\documentclass{article}

\usepackage{amssymb,amsfonts,epsfig}

\title{An ODE--based approach to some Riemann--Hilbert problems
motivated by wave diffraction}

\author{A. V. Shanin}

\newcommand{\ptl}{\partial}

\newtheorem{theorem}{Theorem}

\begin{document}

\maketitle

\begin{abstract}

A novel approach to Riemann--Hilbert problems of particular class is introduced.
The approach is applicable to problems in which the multiplicative jump is set on a half-line.
Such problems are linked to some Wiener--Hopf problems motivated by
diffraction theory.
The new approach is based on
ordinary differential equations: the Riemann--Hilbert problem is reduced to finding a coefficient of
an ordinary differential equation and solving this equation.
The new method leads to an efficient numerical algorithm
and opens a road to new asymptotical and analytical advances.

\end{abstract}

\section{Introduction}

Homogeneous Riemann--Hilbert problems are under consideration. A rigorous formulation of
these problems and some important results related to them can be found in~\cite{Gakhov52}.
These problems are constituted in finding a vector function analytical (or meromorphic) on the complex plane
and having a multiplicative jump on a smooth curve. Traditionally, if the curve coincides with the real axis the problem is referred to as the Wiener--Hopf factorization problem. Wiener--Hopf problems are directly linked with diffraction theory, hydrodynamics and some other branches. Numerous important canonical diffraction problems
(such as diffraction by a half-plane and reflection by waveguide end) can be solved
by the Wiener--Hopf method~\cite{Noble58}. All scalar Wiener--Hopf problems
and only some particular vector Wiener--Hopf problems can be solved
explicitly. The solvable problems include a very important Khrapkov's class \cite{Khrapkov71,Abrahams2007}
admitting commutative factorization. For other problems some approximate
techniques exist (see e.g.~\cite{Abrahams1998}). All problems studied in the papers mentioned above
arise from either wave diffraction theory or solid mechanics.

In the current paper a novel technique is proposed to address the Riemann--Hilbert problems possibly not belonging to known factorizable classes. The approach is based on the ordinary differential equations. While it does not provide
explicit solutions, it opens some new possibilities for asymptotical and analytical studies of the
Riemann-Hilbert problems. Also it leads to an efficient numerical algorithm.

The current paper continues a research started for the problems of wave scattering by a diffraction gratings.
A general framework has been developed for these problems, comprising application of the so-called
embedding formula, spectral equation and the problem of reconstruction of an ODE coefficient.
The current work is dedicated to a possible generalization of the last step of the procedure.

\section{The main result}

\subsection{Problem formulation}

Consider the following matrix Riemann--Hilbert problem.
Let straight contour $C$ connect the point $b$ with $b + i \infty$. The complex
plane is cut along~$C$. The left edge of the cut is denoted by
$C^-$, and the right edge is denoted by~$C^+$.

\begin{figure}[ht]
\centerline{\epsfig{file=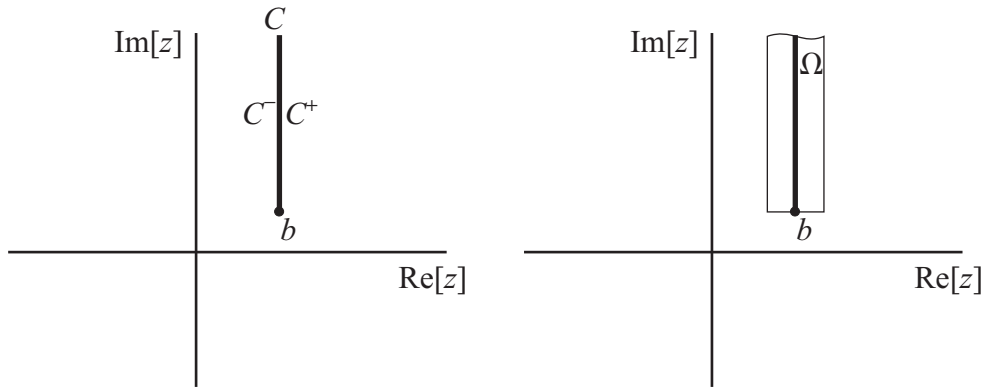}}
\caption{Contour $C$ and domain $\Omega$}
\label{fig01}
\end{figure}

Let there exist an open strip $\Omega$ of non-zero thickness embracing the cut~$C$,
such that the
point $b$ is located on the boundary of $\Omega$ (see Fig.~\ref{fig01}).
Let $2\times 2$ matrix function ${\rm M}(z)$
be analytical in $\Omega$ and continuous at $z = b$.
Let also be
\begin{equation}
{\rm M}(z) \to {\rm I} \quad \mbox{ as } \quad z \to b + i \infty,
\label{eq0101}
\end{equation}
where ${\rm I}$ is a $2 \times 2$ identity matrix.

The Riemann--Hilbert problem studied in this paper is formulated as follows. It is necessary to find
$2 \times 2$ matrix function ${\rm U}$ analytical on the complex plane cut along contour $C$, having
no zeros of the determinant on the cut plane, and obeying the following condition of the cut:
\begin{equation}
{\rm U}^+(z) = {\rm U}^-(z) \, {\rm M} (z) ,
\label{eq0102}
\end{equation}
where ${\rm U}^+(z)$ and ${\rm U}^-(z)$ are the limiting values of the function ${\rm U}$ on the edges
of the cut. Note that the right multiplication by ${\rm M}$ is used in problem formulation.

Additionally, demand that the function  ${\rm U}(z)$ tends to ${\rm I}$ as $|z| \to \infty$.
Thus, a {\em canonical\/} solution is looked for (in Gakhov's sense, \cite{Gakhov52}).
It is assumed also that ${\rm U}$ has no stronger than power growth near $z = b$.

\subsection{Motivation of the problem}

The problem formulation can be considered as a particular case of a standard Riemann--Hilbert problem set on a line going from
$b+ i \infty$ to~$b-i \infty$. Namely, $C$ is a part of this line, and ${\rm M}$ can be defined
on the rest of the line as a unity matrix. However, the author has in mind another link to a standard
matrix factorization problem. Namely, let $b$ belong to the upper half-plane, and let the contour $C'$ coincide
with the real axis.
Let the unknown matrix function ${\rm U}$ be analytical
everywhere except $C'$ and be discontinuous on the
contour $C'$. Let the jump on the contour $C'$ be described by the relation
\begin{equation}
{\rm U}(z+0i) = {\rm U}(z-0i) \, {\rm A} (z) , \qquad z \in C'
\label{eq0102a}
\end{equation}
where ${\rm A}(z)$ is a known matrix coefficient.
Finally, let ${\rm A}$ be an algebraic function of $z$
defined on the whole complex plane $z$ and
having branch points only at $z = \pm b$.
The last condition is quite strong, however this is not a rare occasion in the diffraction--motivated
factorization problems.


Denote by $\bar {\rm U}$ the continuation of ${\rm U}$ from the lower half-plane. According
to (\ref{eq0102a}) this continuation is given by
\begin{equation}
\bar {\rm U}(z) = {\rm U}(z) {\rm A}^{-1} (z), \qquad {\rm Im} [z] >0,
\label{eq0103}
\end{equation}
i.e.\ the function $\bar {\rm U}(z)$ can be continued onto the complex plane cut along~$C$.

Denote by the upper indices $\pm$ the values of $\bar{\rm U}$
and ${\bf A}$ related to the edges of the cut~$C$. Obviously,
\begin{equation}
\bar {\rm U}^+ = {\rm U}(z) ({\rm A}^+ (z) )^{-1},
\qquad
\bar {\rm U}^- = {\rm U}(z) ({\rm A}^- (z) )^{-1} .
\label{eq0104}
\end{equation}
Finally,
\begin{equation}
\bar {\rm U}^+(z) = \bar {\rm U}^-(z)\, {\rm A}^- (z) \, ({\rm A}^+ (z) )^{-1},
\qquad
z \in C
\label{eq0105}
\end{equation}
Thus, function $\bar {\rm U}$ obeys a ``half-line'' Riemann--Hilbert problem of the
type introduced above with
\[
{\rm M}(z) = {\rm A}^- (z)  \, ({\rm A}^+ (z) )^{-1}.
\]
It is possible to reformulate a lot of practical diffraction problems as half--line
Riemann--Hilbert problems.

\subsection{``Ordered exponential'' notations}

Let ${\rm B}(\tau)$ be a $2\times 2$ matrix function analytical
in some domain of complex plane. Define the {\em ordered exponential\/} of ${\rm B}$
(denoted by ${\rm OE}[{\rm B}]$) as a solution of differential equation
\begin{equation}
\frac{\ptl }{\ptl \tau_2} {\rm OE}[{\rm B}](\tau_2 , \tau_1) =
{\rm B}(\tau_2) \cdot {\rm OE}[{\rm B}](\tau_2 , \tau_1)
\label{eq2001}
\end{equation}
with initial condition
\[
 {\rm OE}[{\rm B}](\tau_1 , \tau_1) = {\rm I}.
\]
Here $\tau_1$ and $\tau_2$ belong to the domain
of analiticity of ${\rm B}(\tau)$. Obviously, in this case
$ {\rm OE}[{\rm B}](\tau_2 , \tau_1)$ does not depend on the contour along which
the ordinary differential equation is solved.

If a large area of variation of $\tau_2$ is considered, in which ${\rm B}(\tau_2)$
can have poles, it is necessary to indicate a contour $\gamma$
connecting $\tau_1$ and $\tau_2$, along which equation (\ref{eq2001}) is solved.
In this case the notation looks like ${\rm OE}[{\rm B}](\gamma)$.

Obvious properties of the ordered exponential are as follows:
\begin{equation}
{\rm OE}[{\rm B}](\gamma_1 + \gamma_2) =
{\rm OE}[{\rm B}](\gamma_2) \cdot {\rm OE}[{\rm B}](\gamma_1),
\label{eq2001a}
\end{equation}
\begin{equation}
{\rm OE}[{\rm B}](-\gamma) =
\left( {\rm OE}[{\rm B}](\gamma) \right)^{-1}.
\label{eq2001b}
\end{equation}
Here and below $\gamma_1 + \gamma_2$ denotes concatenation of the contours, $-\gamma$ corresponds
to the contour $\gamma$ passed in the opposite direction.

\subsection{Formulation of the main result of the paper}

\begin{figure}[ht]
\centerline{\epsfig{file=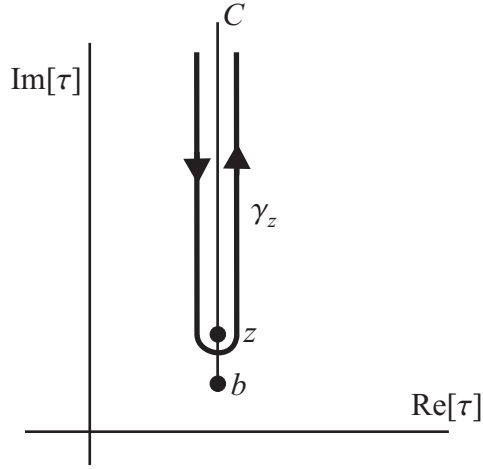}}
\caption{Contours for the OE equation}
\label{fig04}
\end{figure}

\begin{theorem}
Let ${\rm r}(\tau)$ be a function defined and analytical in $\Omega$ and continuous at $\tau = b$
such that for all points $z \in (C \setminus b)$ the following equation is valid:
\begin{equation}
{\rm OE} \left[ \frac{{\rm r}(\tau)}{z-\tau} \right](\gamma_z) = {\rm M}(z) ,
\label{eq2006}
\end{equation}
where contour $\gamma_z$ (shown in Fig.~\ref{fig04}) lies within $\Omega$.
$z$ in the left of (\ref{eq2006}) plays role of a
fixed parameter, while $\tau$ is the independent variable. Then solution of the Riemann--Hilbert
problem formulated above is given by the formula
\begin{equation}
{\rm U}(z) = {\rm OE} \left[
\frac{{\rm r}(\tau)}{z-\tau}
\right](b , b+ i \infty) ,
\label{eq2002}
\end{equation}
where the ordinary differential equation is solved along contour $C$.
\label{th01}
\end{theorem}

{\sc Proof}

\begin{figure}[ht]
\centerline{\epsfig{file=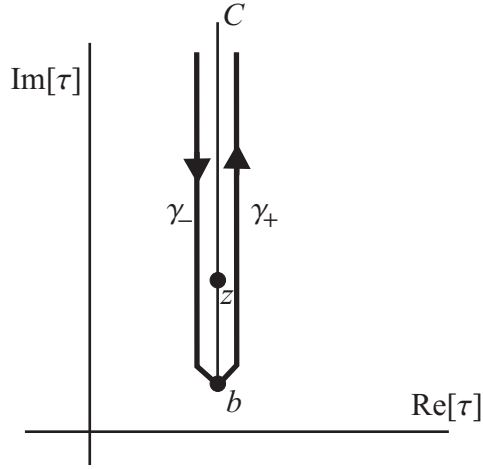}}
\caption{Contours $\gamma_+$ and $\gamma_-$}
\label{fig03}
\end{figure}

Obviously, relation (\ref{eq2002}) defines a function analytical on the complex plane
cut along contour~$C$. Thus, it is necessary to prove (\ref{eq0102}). Let be $z \in C$. Deform contour
$\gamma_z$ as it is shown in Fig.~\ref{fig03}. The contour becomes a concatenation of $\gamma_-$ going from
$+i \infty$ to $b$ along the left edge of the cut, and the contour $\gamma_+$ going from $b$ to $+ i \infty$
along the right edge of the cut.
According to (\ref{eq2001a}), (\ref{eq2001b}), and (\ref{eq2006}),
\begin{equation}
 \left(  {\rm OE} \left[
\frac{{\rm r}(\tau)}{z-\tau}
\right](-\gamma_+) \right)^{-1}
{\rm OE} \left[
\frac{{\rm r}(\tau)}{z-\tau}
\right](\gamma_-) = {\rm M}(z).
\label{eq2010}
\end{equation}
However, according to (\ref{eq2002})
\begin{equation}
{\rm OE} \left[
\frac{{\rm r}(\tau)}{z-\tau}
\right](\gamma_-) = {\rm U}^+(z) ,
\qquad
{\rm OE} \left[
\frac{{\rm r}(\tau)}{z-\tau}
\right](-\gamma_+) = {\rm U}^-(z) ,
\label{eq2011}
\end{equation}
thus (\ref{eq0102}) is proven $\square$

\section{Interpretation of the main result}

Although the proof of Theorem~\ref{th01} is simple and straightforward, it would be useful to
demonstrate another interpretation of this theorem. Let $\beta$ be a point on contour $C$, and let $C_\beta$
be a part of $C$ going from $\beta$ to $+i \infty$. Formulate a family of Riemann--Hilbert problems,
indexed by $\beta$,
having cut $C_\beta$ instead of $C$ (all other conditions of the problem remain the same).
The coefficient of each problem belonging to the family is a restriction of
function ${\rm M}(z)$ on contour~$C_\beta$. Note that ${\rm M}$ does not depend on~$\beta$.

Denote by ${\rm U}(\beta ; z)$ the family of solutions for these problems.
Consider the derivative
\begin{equation}
{\rm W}(\beta ; z) = \frac{\ptl}{\ptl \beta} {\rm U} (\beta ; z).
\label{eq0201a}
\end{equation}
${\rm W} (\beta ; z)$ taken as a function of $z$ is isomonodromical to ${\rm U} (\beta ; z)$
in the following sense: it is analytical on the plane with the same cut $C_\beta$, and the multiplicative jump on this cut is described
by the same coefficient ${\rm M}(z)$. Moreover, for each $\beta$ function ${\rm W}(\beta ; z)$ decays as
$|z| \to \infty$, and this function grows no stronger than a power function near $z = \beta$.

The isomonodromy between ${\rm W}$ and ${\rm U}$ can be used as follows. One can prove that
\begin{equation}
{\rm W}(\beta ; z) = {\rm R} (\beta ; z) \, {\rm U} (\beta)
\label{eq0202}
\end{equation}
where ${\rm R}(\beta ; z)$ is a rational function of $z$ for each fixed $\beta$.
To prove this, define ${\rm R}$ as
\[
{\rm R} = {\rm W} \, {\rm U}^{-1}.
\]
Obviously, ${\rm R}$ has no singularities on the plane cut along $C_\beta$, except maybe the point $z = \beta$.
The discontinuity of ${\rm W}$ on $C_\beta$ is described by the relation
\[
{\rm R}^+ (\beta ; z) = {\rm W}^+ (\beta ; z) \, ({\rm U}^+(\beta ; z))^{-1} =
{\rm W}^- (\beta ; z) \, {\rm M}(z) \, {\rm M}^{-1} (z) ({\rm U}^-(\beta ; z))^{-1}
=
\]
\begin{equation}
{\rm W}^- (\beta ; z) \, ({\rm U}^-(\beta ; z))^{-1}
= {\rm R}^- (\beta ; z),
\label{eq0203}
\end{equation}
i.e.\ ${\rm R}$ is continuous on $C_\beta$. Thus, ${\rm R}$ is single valued, has a power singularity
at $z = \beta$, and decays as $|z| \to \infty$. Therefore it is rational.
A detailed study of the singularity at $z = \beta$ shows that ${\rm R}$ should have form
\begin{equation}
{\rm R} (\beta ; z) = \frac{{\rm r} (\beta)}{z-\beta} .
\label{eq0204}
\end{equation}
Taking into account that the Riemann--Hilbert problem degenerates as $\beta \to b+ i \infty$, and
${\rm U}(b+i\infty ; z) = {\rm I}$, the representation (\ref{eq2002}) follows from
(\ref{eq0202}) rewritten as a differential equation
\begin{equation}
\frac{\ptl}{\ptl \beta}{\rm U}(\beta ; z) = {\rm R} (\beta ; z) \, {\rm U} (\beta ; z).
\label{eq0205}
\end{equation}

{\bf Note 1.}
The method proposed here is based on embedding of the considered Riemann--Hilbert problem
into a family having parameter $\beta$ such that the differentiation with respect to the parameter
provides an isomonodromy. In this context an elementary problem can be considered, for which ${\rm M}$ is
a constant matrix with respect to~$z$ (or there are several cuts $C_j$ and different
constant matrices ${\rm M}_j$ defined on them).
It is not necessary to embed this problem into a family
of other problems, since differentiation $d / d z$ provides an isomonodromy itself. The argument similar to
the one provided above shows that ${\rm U}$ in this case is a solution of a Fuchsian equation,
and ${\rm M}_j$ are the monodromy matrices of this equation.

{\bf Note 2.} Variation of the end point of the cut is not the only possible way to organize a family, in which
the differentiation with respect to a parameter provides an isomonodromy. Another way is to study matrix
${\rm M}(\alpha; z)$ depending on some parameter $\alpha$ and construct a differential operator
\[
H \equiv a_1 (\alpha ; z) \frac{\ptl}{\ptl \alpha} + a_2(\alpha ; z) \frac{\ptl}{\ptl z}
\]
such that $H[{\rm M}] = 0$. Then it is possible prove that
\[
H[{\rm U}](\alpha ; z) = {\rm R}(\alpha ; z) \, {\rm U}(\alpha ; z),
\]
where ${\rm R}$ is a single-valued function with respect to $z$.

\section{Examples}
\subsection{Scalar example}

Consider a scalar example.
Let it be necessary to
find function $u(z)$ obeying the following restrictions:
\begin{itemize}

\item
The function should be regular on the complex plane cut along the
line $ C = (b , b + i \infty)$.

\item
The values on the edges of the cut should be connected via the relation
\begin{equation}
u^+(z) = u^-(z)\,m(z),
\label{eq0101a}
\end{equation}
where $u^+$ are the values on the right edge of the cut, $u^-$ are the values on the left edge of the cut, $m(z)$, $z \in \Omega$ is known function tending to 1 at infinity.

\item
The function should grow no faster than algebraically at $b$ and tend to 1 at infinity.

\end{itemize}


Obtain a ``traditional'' solution for this problem. Note that a multiplicative jump
of $u$ corresponds to an additive jump of $\log u$. Thus, the solution
can be obtained by using Cauchy's integral:
\begin{equation}
\log u(z) =
- \frac{1}{2\pi i} \int \limits_b^{b+i \infty}  \frac{ \log ( m(\xi) )}{\xi-z} d\xi.
\label{eq0201}
\end{equation}

Now construct a ``new'' solution for this problem. According to a scalar analog of Theorem~\ref{th01},
it is necessary to find a function $r(\tau)$ defined in $\Omega$, such that an analog of
(\ref{eq2006}) is fulfilled:
\begin{equation}
{\rm OE} \left[ \frac{r(\tau)}{z-\tau} \right](\gamma_z) = m(z) ,
\label{eq0601}
\end{equation}
This equation can be simplified using the fact that
\begin{equation}
{\rm OE}  \left[ \frac{r(\tau)}{z-\tau} \right](\gamma_z) =
\exp\left\{ \int_{\gamma_z} \frac{r(\tau)}{z-\tau} d \tau \right\}.
\label{eq0602}
\end{equation}
Thus one can find $r(\tau)$ explicitly:
\begin{equation}
r(\tau) = - \frac{\log (m(\tau))}{2\pi i}.
\label{eq0404}
\end{equation}
Using the scalar analog of (\ref{eq2002})
\begin{equation}
u(z) =
\exp\left\{-
\int \limits^{b+i \infty}_b
\frac{ r(\tau)}{z-\tau} d\tau
\right\},
\label{eq0602a}
\end{equation}
get the solution coinciding with (\ref{eq0201}).

\subsection{Chebotarev--Khrapkov matrices}

The most important class of matrices admitting explicit factorization is the class
introduced by Chebotarev \cite{Chebotarev56} and Khrapkov \cite{Khrapkov71}. Later on, these matrices have been studied by numerous researchers. Among the most important works there are \cite{Rawlins81,Daniele84,Jones84,Antipov02}.

Here the Khrapkov's $2 \times 2$ matrices are taken as an example. These matrices have the following form:
\begin{equation}
{\rm M}(z) = c(z) {\rm I} + p(z) {\rm L} , \qquad
{\rm L}(z) \equiv \left(   \begin{array}{cc}
l(z) & m(z) \\
n(z) & -l(z)
\end{array} \right)
\label{eq0701}
\end{equation}
where $l(z)$, $m(z)$, $n(z)$ are some polynomials such that
\begin{equation}
f(z) \equiv l^2(z) + m(z) n(z)
\label{eq0702}
\end{equation}
is a polynomial of degree not higher than~2.
It is assumed that ${\rm M} \to {\rm I}$ as $z \to b+ i \infty$.

It is known that the Khrapkov's matrices can be factorized commutatively, but for the price of more flexible behaviour at infinity. Thus, the coefficient
\[
{\rm R}(\tau ; z) = \frac{{\rm r}(\tau)}{z-\tau}
\]
in (\ref{eq2002}) and (\ref{eq2006})
should be replaced by a form not necessarily decaying as $|z| \to \infty$:
\begin{equation}
{\rm R}(\tau ; z) =  \frac{\xi(\tau)}{z - \tau} {\rm I} + \frac{\eta(\tau)}{z-\tau} {\rm L}(z),
\label{eq0703}
\end{equation}
where $\xi$ and $\eta$ are some scalar functions.

The advantage of using the Khrapkov's matrices is their commutativity which leads
to the closed-form formula for the ordered exponential:
\begin{equation}
{\rm OE}[{\rm R}(\tau ; z)](\gamma) = \exp\{ \hat \xi \}
\left(
{\rm cosh}\left( \sqrt{f(z)} \, \hat \eta \right) {\rm I} +
{\rm sinh}\left( \sqrt{f(z)} \, \hat \eta \right)
\frac{ {\rm L}(z) }{ \sqrt{f(z)} }
\right),
\label{eq0704}
\end{equation}
where
\begin{equation}
\hat \xi  = \int_\gamma \frac{\xi(\tau)}{z-\tau} d\tau,
\qquad
\hat \eta = \int_\gamma \frac{\eta(\tau)}{z-\tau} d\tau,
\label{eq0705}
\end{equation}

Thus, in the analog of (\ref{eq2006})
\begin{equation}
{\rm OE} [ {\rm R}(\tau ; z) ] (\gamma_z) = {\rm M} (z)
\label{eq0706}
\end{equation}
contour $\gamma_z$ can be changed to a small closed loop encircling
the point $z$. According to this, (\ref{eq0706}) can be solved explicitly:
\begin{equation}
\xi(\tau) = \frac{i}{4 \pi} \log (c^2(\tau) - f(\tau) p^2(\tau)),
\label{eq0707}
\end{equation}
\begin{equation}
\eta(\tau) = \frac{i}{4 \pi \sqrt{f (\tau)}} \log \left(
\frac{c(\tau) + p(\tau) \sqrt{f (\tau)}}{c(\tau) - p(\tau) \sqrt{f (\tau)}}
\right).
\label{eq0708}
\end{equation}
Finally, the solution of the problem is as follows:
\begin{equation}
{\rm U}(z) = \exp\{ \bar \xi \}
\left(
{\rm cosh}\left( \sqrt{f(z)} \, \bar \eta \right) {\rm I} +
{\rm sinh}\left( \sqrt{f(z)} \, \bar \eta \right)
\frac{ {\rm L}(z) }{ \sqrt{f(z)} }
\right),
\label{eq0709}
\end{equation}
where
\begin{equation}
\bar \xi  = -\int \limits_b^{b + i \infty} \frac{\xi(\tau)}{z-\tau} d\tau,
\qquad
\bar \eta = -\int \limits_b^{b + i \infty} \frac{\eta(\tau)}{z-\tau} d\tau.
\label{eq0710}
\end{equation}
Solution in the form (\ref{eq0709}) coincides with the one described in \cite{Khrapkov71}.

\subsection{A link to Weinstein--class problems}

\begin{figure}[ht]
\centerline{ \epsfig{file=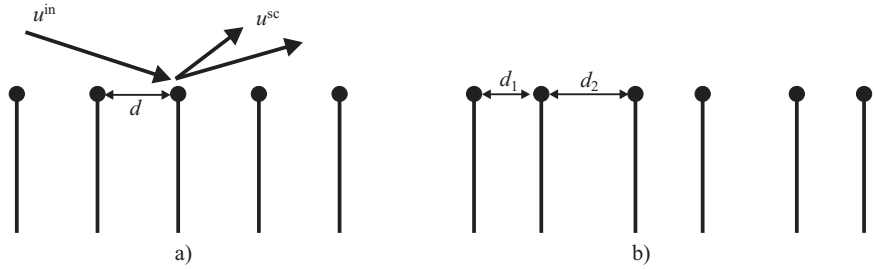, width = 11.5cm} }
\caption{Geometry of Weinstein--class problems}
\label{fig08}
\end{figure}

The current paper can be considered as a generalization of a previous work by author dedicated to
Weinstein--class problems. Problems of scattering by periodic diffraction gratings composed of
absorbing screens have been studied. The first problem having the geometry shown in Fig.~\ref{fig08}~a) is a classical
Weinstein's problem \cite{Weinstein0}. The period of the system consists of a single screen.
The incident wave has wavenumber $k_0$ (such that $k_0 d \gg 1$) and a grazing incidence angle.
The problem is considered in the parabolic approximation.
In \cite{Weinstein1} this problem has been reformulated as the following scalar OE equation:
\begin{equation}
{\rm OE}[r (\tau)/ \sqrt{z-\tau}](\gamma_z) = 1 - e^{i z},
\label{eq2301}
\end{equation}
where $\gamma_z$ is the contour introduced above. This equation can be reduced to an integral equation
with difference kernel and then solved explicitly. It is worth to note that this problem has been solved
in \cite{Weinstein0} by using Wiener--Hopf method.

In the second problem of Weinstein's class (see Fig.~\ref{fig08}~b)) the period of the grating contains
two screens. In \cite{Weinstein2} this problem has been reduced to matrix OE equation
\begin{equation}
{\rm OE}\left[
\frac{{\rm r}(\tau)}{\sqrt{z-\tau}}
\right] (\gamma_z) = {\rm M}(z),
\label{eq2302}
\end{equation}
\begin{equation}
{\rm M}(z) = \left(  \begin{array}{cc}
1 & -\exp \{ i z d_2 / (d_1 + d_2)  \} \\
- \exp \{ i z d_1 / (d_1 + d_2) \} & 1
\end{array} \right).
\label{eq2303}
\end{equation}
The same diffraction problem can be reduced to
a Wiener--Hopf problem in which the matrix
\begin{equation}
{\rm K}(\beta ; y ) = {\rm M} (\beta - y^2 (d_1 + d_2)/ (2 k_0))
\label{eq2304}
\end{equation}
should be factorized with respect to variable $y$ for any real value of parameter~$\beta$.

Application of the equation (\ref{eq2302}) enabled one to construct an efficient numerical
solution of the diffraction problem having geometry shown in  Fig.~\ref{fig08}~b. Moreover,
using the asymptotic analysis of equation (\ref{eq2302}) it became possible to prove that the
reflection coefficient for this structure tends to $-1$ as the angle on incidence tends to zero.
This result is practically important since it guarantees high values of Q-factor in
corresponding resonator problems~\cite{Shabalina10}.

\section{Solving (\ref{eq2006}) numerically}

\subsection{Parametrization of equation (\ref{eq2006})}

Here equation (\ref{eq2006}) is addressed.
It is shown below that the unknown matrix ${\rm r} (\tau)$ depends on two
unknown scalar parameters. These parameters are revealed here, and the problem for finding these
parameters is formulated in the simplest.

Let the logarithm of the eigenvalues of ${\rm M}$ (and the eigenvalues themselves)
be distinct everywhere on~$C$. Represent ${\rm M}$ in the form
\begin{equation}
{\rm M}(\tau) = {\rm T}(\tau) \cdot  \left( \begin{array}{cc}
\lambda_1 (\tau) & 0 \\
0 & \lambda_2 (\tau)
\end{array}  \right) \cdot {\rm T}^{-1} (\tau)
\label{eq0801}
\end{equation}
Here ${\rm T}$ is the matrix whose columns are eigenvectors of ${\rm M}$. Matrix ${\rm T}$ depends on
two scalar parameters. Almost everywhere it can be parametrized as follows:
\begin{equation}
{\rm T}(\tau)  = \left( \begin{array}{cc}
1          & 1 \\
t_1 (\tau) & t_2 (\tau)
\end{array} \right)
\label{eq0802}
\end{equation}
Values $\lambda_{1,2}$ and $t_{1,2} (\tau)$ are known.

The left-hand side of (\ref{eq2006}) can be rewritten as follows:
\begin{equation}
{\rm OE}\left[ \frac{{\rm r}(\tau)}{z-\tau} \right](\gamma_z) =
{\rm F} \cdot {\rm OE}\left[ \frac{{\rm r}(\tau)}{z-\tau} \right](\sigma_z)
\cdot {\rm F}^{-1}
\label{eq0802a}
\end{equation}
where $\sigma_z$ is a loop of a small radius $\epsilon$ starting and ending
at $z + i \epsilon$ and encircling $z$ in the positive direction, and
\[
{\rm F}
= {\rm OE}\left[ \frac{{\rm r}(\tau)}{z-\tau} \right](z + i \infty , z+i \epsilon).
\]
Let ${\rm r}(\tau)$ be represented in the form
\begin{equation}
{\rm r}(\tau) = {\rm H}(\tau) \, \left(   \begin{array}{cc}
\zeta_1(\tau) & 0 \\
0 & \zeta_2 (\tau)
\end{array} \right)
\,
{\rm H}^{-1} (\tau)
\label{eq0803}
\end{equation}
As $\epsilon \to 0$
\begin{equation}
{\rm OE}\left[ \frac{{\rm r}(\tau)}{z-\tau} \right](\sigma_z)
\rightarrow
{\rm H}(z) \cdot
\exp \left\{
- 2\pi i
\left(   \begin{array}{cc}
\zeta_1(z) & 0 \\
0 & \zeta_2 (z)
\end{array} \right) \right\}
\cdot
{\rm H}^{-1}(z).
\label{eq0804}
\end{equation}

It follows from (\ref{eq0802}) that the eigenvalues of ${\rm r}(z)$ are defined
by the eigenvalues of ${\rm M}(z)$:
\begin{equation}
\zeta_{1,2}(z) = \frac{i}{2\pi} \log(\lambda_{1,2}(z)).
\label{eq0805}
\end{equation}
Thus, to find ${\rm r}$ one needs only to find ${\rm H}(\tau)$.

Again, the columns of ${\rm H}(\tau)$ are the eigenvectors of ${\rm r}$. Almost everywhere
matrix ${\rm H}$ can be parametrized as follows:
\begin{equation}
{\rm H}(\tau)  = \left( \begin{array}{cc}
1          & 1 \\
h_1 (\tau) & h_2 (\tau)
\end{array} \right)
\label{eq0806}
\end{equation}
So our aim is to find $h_{1,2} (\tau)$.

\begin{figure}[ht]
\centerline{ \epsfig{file=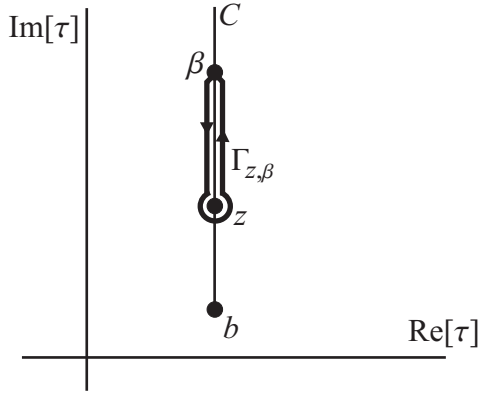} }
\caption{Contour $\Gamma_{z,\beta}$}
\label{fig20}
\end{figure}

For $z, \beta \in C$ , ${\rm Im}[\beta] > {\rm Im}[z]$ define the function
\begin{equation}
{\rm V}(\beta ; z) =
{\rm OE}\left[ \frac{{\rm r}(\tau)}{z-\tau} \right](\Gamma_{z,\beta})
\label{eq0807}
\end{equation}
where contour $\Gamma_{z, \beta}$ is shown in Fig.~\ref{fig20}.
Obviously, ${\rm V}$ obeys the equation
\begin{equation}
\frac{\ptl}{\ptl \beta}  {\rm V}(\beta  ; z) =
\frac{1}{z-\beta}
\left(
{\rm r}(\beta) \, {\rm V}(\beta ; z) -
{\rm V}(\beta ; z) \, {\rm r}(\beta)
\right) .
\label{eq0808}
\end{equation}

Since ${\rm V}$ is adjoint to to
\[
{\rm OE}\left[ \frac{{\rm r}(\tau)}{z-\tau} \right](\sigma_z),
\]
its eigenvalues are the same, and it can be written as follows:
\begin{equation}
{\rm V}(\beta ; z) =
{\rm Q}(\beta ; z) \cdot
\left( \begin{array}{cc}
\lambda_1 (z) & 0 \\
0 & \lambda_2 (z)
\end{array}  \right)
\cdot
({\rm Q}(\beta ; z))(z)^{-1}.
\label{0809}
\end{equation}
where
\begin{equation}
{\rm Q}(\beta ; z)  = \left( \begin{array}{cc}
1          & 1 \\
q_1 (\beta ; z) & q_2 (\beta ; z)
\end{array} \right).
\label{eq0810}
\end{equation}
Taking $\beta \to z$ obtain the relation ${\rm Q}(z ; z) = {\rm H}(z)$, and thus
\begin{equation}
q_{1,2} (z;z) = h_{1,2}(z).
\label{eq0811}
\end{equation}
Taking $\beta \to b + i \infty$ obtain ${\rm Q}(\infty ; z) = {\rm M}(z)$, and thus
\begin{equation}
q_{1,2} (b + i \infty ; z) = t_{1,2} (z).
\label{eq0812}
\end{equation}

Elementary calculations demonstrate that equation (\ref{eq0808}) is equivalent to the following system of
differential equations:
\begin{equation}
\frac{\ptl }{\ptl \beta} q_{1,2} (\beta , z) =
- \frac{\left(\zeta_1(\beta) - \zeta_2(\beta) \right)
\left( q_{1,2}(\beta ; z) -  h_1(\beta) \right)
\left( q_{1,2}(\beta ; z) -  h_2(\beta) \right)
}{(z - \beta) \left(h_1(\beta) - h_2 (\beta)\right)}
\label{eq0813}
\end{equation}
Note that these equations are independent from each other. Note also that they are Riccati equations. 
Finally, the reduced problem can be formulated as follows:

\begin{em}
For given functions $\zeta_{1,2} (\tau)$, $t_{1,2} (\tau)$ , $\tau \in C$
find functions $h_{1,2}(\tau)$, $\tau \in C$ such that for any $z \in C$ equation (\ref{eq0813})
possesses a solution $q_{1,2}(\beta, z)$ on the segment $\beta \in (z , b+ i \infty)$ of $C$ obeying boundary
conditions (\ref{eq0811}), (\ref{eq0812}).
\end{em}

\subsection{Numerical algorithm}

The process of solving of the Riemann--Hilbert problem using the new method consists of two major
stages. On the first stage the coefficient ${\rm r}(\tau)$ is found on the contour $C$
(it is shown above that it is sufficient to find two scalar functions $h_{1,2} (\tau)$).
The coefficient is found by solving equation (\ref{eq2006}) or (\ref{eq0813}).
On the second stage equation (\ref{eq2002}) is solved to find function ${\rm U}(z)$ at some set of points
$z_k$. The points at which the ${\rm U}$ should be computed do not necessarily belong to $C$.
Moreover, in diffraction applications $C$ is a half-line parallel to the imaginary axis, while it is
necessary to find ${\rm U}$ on the real axis.

\begin{figure}[ht]
\centerline{ \epsfig{file=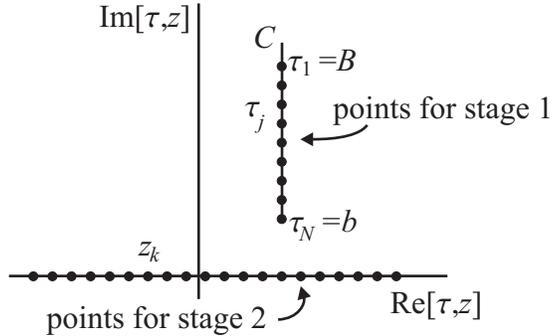} }
\caption{Nodes for numerical computations}
\label{fig21}
\end{figure}

For the first stage contour $C$ is meshed, i.e.\ a dense array of nodes $\tau_j$ is placed on it. Since the contour is infinite, it is necessary to choose a reasonable upper limit $B$ instead of $b + i \infty$
(see Fig.~\ref{fig21}). The nodes are numbered by indices $j = 1\dots N$, where $1$ corresponds to the node
$\tau_1 = B$ (infinity), and index $N$ corresponds to $\tau_N = b$.

At each $\tau_j$ matrix ${\rm M}(\tau)$ is represented in the form (\ref{eq0801}), (\ref{eq0802}), i.e.
the values $\lambda_{1,2}(\tau_j)$ and $t_{1,2}(\tau_j)$ are computed. The values $\zeta_{1,2}(\tau_j)$
are then computed by applying formula (\ref{eq0805}). The aim is to compute the values $h_{1,2}(\tau_j)$.

For the ``infinity'' point $\tau_1$ the following values are assigned:
\begin{equation}
h_{1,2}(\tau_1) = t_{1,2}(\tau_1).
\label{0901}
\end{equation}
This is a natural choice for the asymptotics of the unknown coefficient, since
${\rm M}(\tau) \to {\bf I}$ as $\tau \to b + i \infty$.

Then a loop over $j = 2 \dots N$ is performed. At the $j$th
step of the loop  the values $h_{1,2}(\tau_j)$ should be computed. Thus, when the
$j$th step is performed all values $h_{1,2}(\tau_1) \dots h_{1,2}(\tau_{j-1})$ are
already found. On the $j$th step of the loop equations (\ref{eq0813})
are solved on the segment $(\tau_1 , \tau_{j-1})$ for $q_{1,2}$
using Runge--Kutta~4 method. The initial values are set as
\begin{equation}
q_{1,2}(\tau_1) = t_{1,2}(\tau_j)
\label{eq0902}
\end{equation}
following from (\ref{eq0812}). The values $q_{1,2}(\tau_{j-1})$
are found. Then equations (\ref{eq0813}) are solved on the segment
$(\tau_{j-1}, \tau_j)$ (only one step is performed) by Euler
method. This method does not require the values of the right-hand side at the end of the
segment to be known. The values $q_{1,2}(\tau_j)$ are found. The assignment
\begin{equation}
h_{1,2}(\tau_j) = q_{1,2}(\tau_j)
\label{eq0903}
\end{equation}
is performed, following from (\ref{eq0811}).

On the second stage for each $z_k$ the array
${\rm r}(\tau_j)/(z_k-\tau_j)$ is computed by (\ref{eq0803}), (\ref{eq0806}).
Then the value
\[
{\rm U}(z_k) =
{\rm OE}\left[ \frac{{\rm r}(\tau)}{z_k - \tau}  \right] (\tau_N , \tau_1)
\]
is is found by Runge--Kutta~4 method.

To demonstrate the feasibility of the new method the Khrapkov's matrix of the simplest
form has been taken:
\begin{equation}
{\rm M}(z) = {\bf I} + \frac{1}{z^2}
\left( \begin{array}{cc}
1 & z \\
-z & -1
\end{array} \right)
\label{eq0904}
\end{equation}
The matrix is set on contour $C = (2i , i\infty)$.
The aim was to find the unknown function ${\rm U}$ at the points $z_k$ belonging to the line
${\rm Im}[z] = 1.8$.

The results obtained by the new method can be compared with the exact solution
provided by Khrapkov's formulae (\ref{eq0709}), (\ref{eq0710}). One has to take into account, however,
that Khrapkov's formulae lead to the solution ${\rm U}_{\rm kh}(z)$ that tends to some constant
${\rm U_\infty} \ne {\rm I}$ as $|z| \to \infty$. Thus, the Khrapkov's solution should be corrected as follows:
\begin{equation}
{\rm U} = ({\rm U}_{\infty})^{-1} ({\rm U}_{\rm kh})(z)
\label{eq0905}
\end{equation}

\begin{figure}[ht]
\centerline{ \epsfig{file=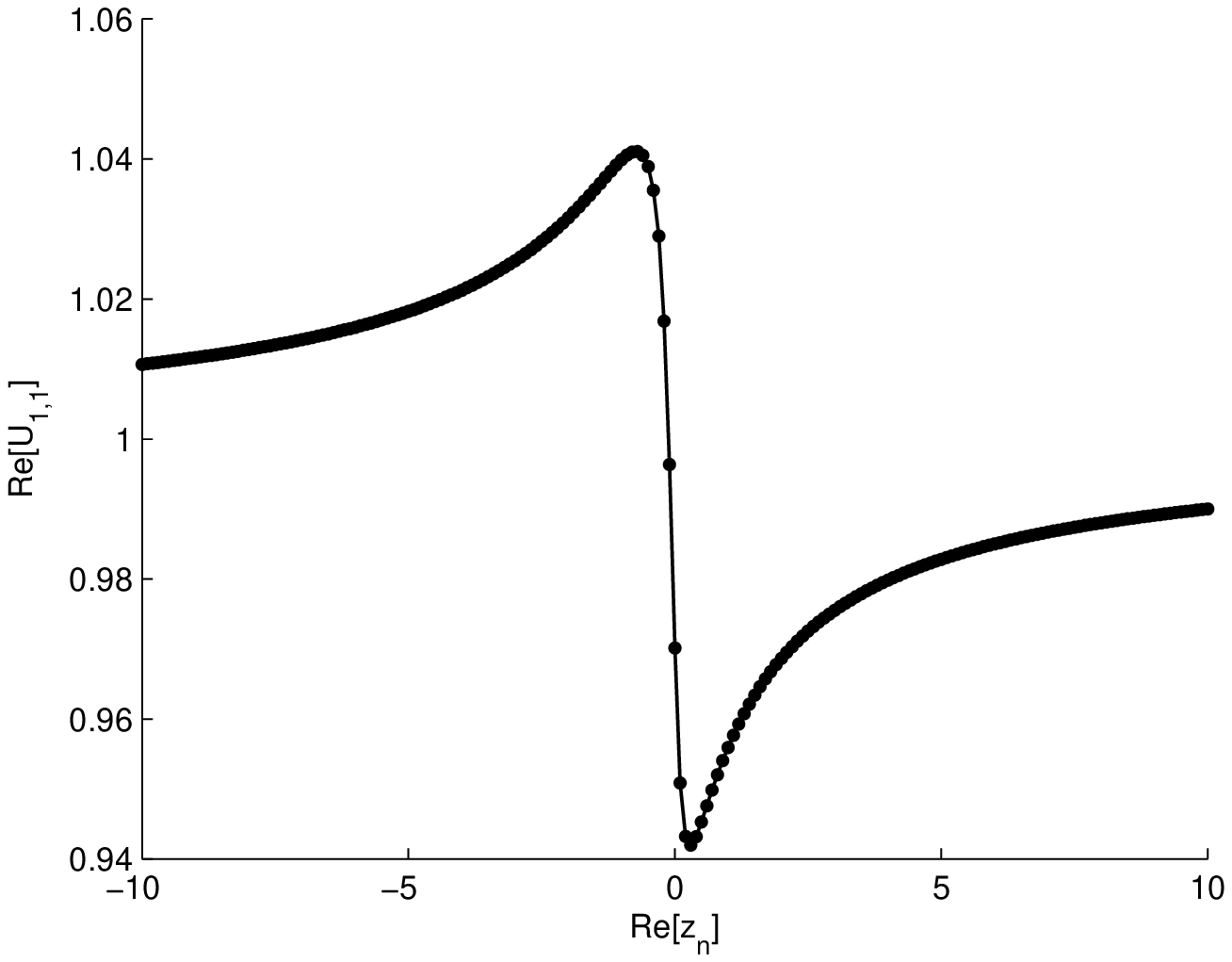,width=6cm},\epsfig{file=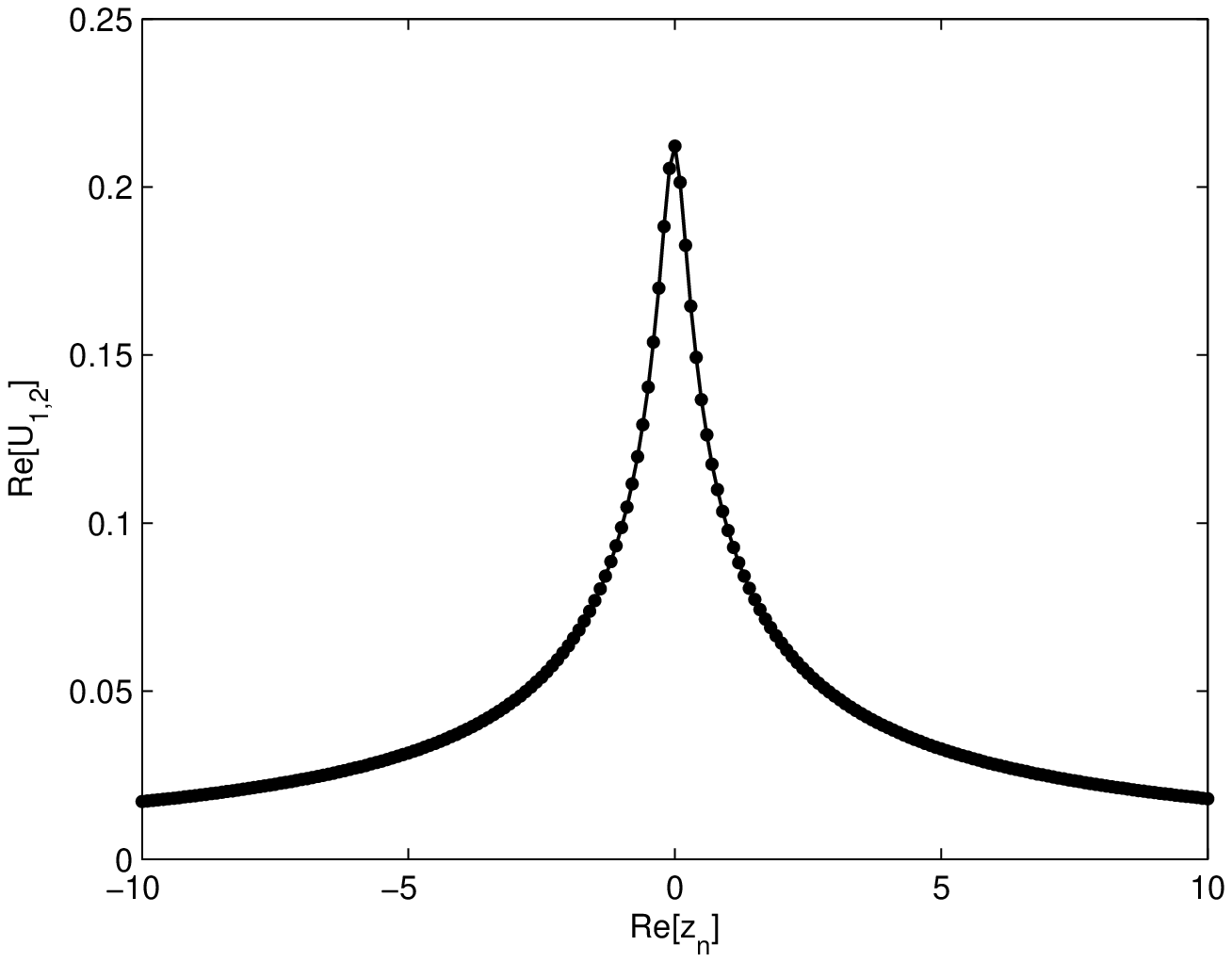,width=6cm}  }
\caption{Results of numerical computations: ${\rm Re[U]}_{1,1}$ (left),
${\rm Re[U]}_{1,2}$ (right).}
\label{fig22}
\end{figure}

For the computations the value of $B$ has been chosen to $80i$, the step on the contour $C$
was equal to $0.02$. The results for ${\rm Re [U]}_{1,1}$ and ${\rm Re [U]}_{1,2}$
are shown in Fig.~\ref{fig22}. Solid lines represent the Khrapkov's solution corrected by
formula (\ref{eq0905}), dots represent the solution obtained by applying the novel ODE-based method.
One can see that the agreement between two exact formula and the new method is very nice.

\section{Conclusion}

The paper presents a new analytical result concerning Riemann--Hilbert problem with the
multiplicative jump set on a half-line. The solution can be represented
in the form (\ref{eq2002}) with the unknown coefficient ${\rm r}(\tau)$ obeying equation
(\ref{eq2006}). Unfortunately, the new result does not lead to explicit solving of a general matrix
factorization problem, but it provides an efficient numerical technique and gives some
hopes to achieve further progress in asymptotical and analytical studies. The efficiency
and robustness of the numerical method is linked with the fact that the equation (\ref{eq2006})
has the ``Volterra'' structure, i.e.\ the unknown coefficient can be found point-by-point. The asymptotical advantages of the new method has already been demonstrated for similar problems of Weinstein type. Finally, some analytical progress, hopefully,
can be achieved by studying the analytical structure of the auxiliary function ${\rm r}(\tau)$ for the case of rational or algebraic matrix ${\rm M}$.

The new method brings the Riemann--Hilbert problem
back into the context of the ordinary differential equations (from which it initially emerged),
but in a slightly unusual form. Instead of studying the derivative with respect to the main variable
$z$ the derivative with respect to the parameter $\beta$ is considered.

A new method is proposed for Riemann--Hilbert problems belonging to a
particular class.
As it is shown in the paper, the class naturally arises as a reformulation of a problem
of factorization of a matrix ${\rm A}$ on the real axis if the matrix is
given as an algebraic function having one branch point in the upper half-plane (or in the lower half-plane). Such case is quite common in the diffraction theory, but it restricts the class of
problems to those having only one propagation speed. If there exist several wave modes
with different speeds the new methods needs to be generalized to admit several half-infinite
contours~$C_j$. Such a generalization can be obtained by allowing the coefficient ${\rm R}$
to have several pole singularities.

The paper considers only matrices $2\times 2$, however most of the results can be easily
converted for the case of matrix of arbitrary dimension.

\section{Acknowledgements}

Author is grateful to Prof.~R.V.~Craster and Prof.~A.Ya.~Kazakov for valuable discussions.
The work has been supported by Russian Government grant 11.G34.31.0066,
Scientific school grant 2631.2012.2 and the Russian Foundation for Basic Research
grant 12-02-00114.

\end{document}